\documentclass[11pt]{amsart}
\usepackage{latexsym}
\usepackage{amsmath}
\usepackage{url}

\newcommand{\PD}{\partial}

\newcommand{\wt}{\widetilde}
\newcommand{\Ac}{\mathcal{A}}

\newcommand{\Cc}{\mathcal{C}}

\newcommand{\Rb}{\mathbb{R}}

\newcommand{\R}{\rangle}
\newcommand{\Beq}{\begin{equation}}
\newcommand{\Eeq}{\end{equation}}
\newcommand{\beq}{\begin{equation*}}
\newcommand{\eeq}{\end{equation*}}
\newcommand{\bal}{\begin{align}}
\newcommand{\eal}{\end{align}}

\newcommand{\g}{\gamma}

\newcommand{\n}{\nabla}

\newcommand{\bp}{\begin{prob}}
\newcommand{\ep}{\end{prob}}
\newcommand{\bpr}{\begin{proof}}
\newcommand{\epr}{\end{proof}}

\usepackage{amssymb,amsmath,amsfonts,amsthm}
\newtheorem{thm}{Theorem}

\newtheorem{prop}{Proposition}

\theoremstyle{definition}

\newtheorem{prob}{Problem}

\begin{document}
\author{V.~Krishnan}
\title{A Support Theorem for the Geodesic Ray Transform of Functions}
\footnotetext{\textit{Math Subject Classifications.}
47G10, 47G30, 53B21.}
\footnotetext{\textit{Keywords and Phrases.}
Geodesic ray transform, Support theorems, Microlocal analysis.}
\begin{abstract}
Let $(M,g)$ be a simple Riemannian manifold. Under the assumption that the metric $g$ is real-analytic, it is shown that if the geodesic ray transform of a function $f\in L^{2}(M)$ vanishes on an appropriate open set of geodesics, then $f=0$ on the set of points lying on these geodesics. The approach is based on analytic microlocal analysis.
\end{abstract}
\maketitle
\section{Introduction}
Let $(M,\PD M, g)$ be a smooth $n$-dimensional compact Riemannian manifold with boundary. In this paper we study the geodesic ray transform of functions. This is defined as follows: Let $\g:[0,l(\g)]\to M$ be a geodesic joining boundary points, where $l(\g)$ is the length of this geodesic. The geodesic ray transform $I$ of a function $f$ along $\g$ is defined as 
\[ If(\g)=\int\limits_{0}^{l(\g)}f(\g(t))dt,\]

We address the question whether the geodesic ray transform of a function $f\in L^{2}(M)$ over an appropriate open set $\Ac$ of geodesics uniquely determines the function $f$ on the set of points lying on the geodesics of $\Ac$. The results of this kind are known as support theorems in integral geometry literature. 

One needs additional restrictions on the metric even to prove injectivity results for this transform as the following counterexample shows \cite{SU3}: Consider the unit sphere with a small disk excised out from its east pole making the resulting manifold a smooth manifold with boundary. Now consider a function $f=1$ on a small disk centered at the north pole and $f=-1$ on a symmetrical disk centered at the south pole. Then the geodesic ray transform vanishes identically, but $f$ is not $0$.

One assumption is to assume that the Riemannian manifold $(M,\PD M, g)$ is simple.
A compact Riemannian manifold with boundary $(M,\PD M,g)$ is simple if
\begin{enumerate}
\item[(a)] The boundary $\PD M$ is strictly convex: $\langle \n_{\xi}\nu,\xi\R > 0$ for $\xi \in T_{x}(\PD M)$. Here $\nu$ is the unit outward normal to the boundary.
\item[(b)] The map $\exp_{x}:\exp_{x}^{-1}M\to M$ is a diffeomorphism for each $x\in M$.
\end{enumerate}

It is well known that simple Riemannian manifolds are diffeomorphic to closed balls in $\Rb^{n}$ \cite{Sh1}. So from now on, without loss of generality, we will assume that our manifold $M$ is a closed ball in $\Rb^{n}$.

It is known from the work of Mukhometov that on a simple Riemannian manifold, the geodesic ray transform is injective. We refer to Sharafutdinov's book \cite{Sh1} for this and other results related to the geodesic ray transform. However a support theorem for this transform has been an open question. In this paper, under the additional assumption of real-analyticity of the simple Riemannian metric, we prove a support theorem. We use analytic microlocal analysis to prove our results. The tools of microlocal analysis were introduced in the Radon transform setting by V. Guillemin. See \cite{GS} for more details.  Based on this work, Boman and Quinto in \cite{BQ1} proved a support theorem for the generalized Radon transform integrated against real analytic weights by using a microlocal version of unique continuation of analytic functions due to Kawai-Kashiwara-H\"ormander. Our support theorem relies on this unique continuation result.

There are several support theorems in integral geometry based on the application of the aforementioned result of Kawai-Kashiwara-H\"ormander. We refer the reader to the following papers for some of the references \cite{B1, B2, BQ1, BQ2, GQ1, GQ2, Q1, Q2, QZ}. The geometric considerations and the transform that we study in this paper makes our result different from these earlier known support theorems.

We make a remark before stating our main result. If $(M,g)$ is a simple Riemannian manifold with $g$ real-analytic, we can extend the manifold slightly to a larger manifold $\wt{M}$ with $M$ contained in the interior of $\wt{M}$, and the metric $g$ real-analytically to $\wt{M}$ (we call the extended metric $g$) such that $(\wt{M},g)$ is simple. From now on, we will assume that $\wt{M}$ is a simple manifold as above. Also when we talk about geodesics, we will assume that they are maximal geodesics joining boundary points.

We now state our main result.
\begin{thm}\label{S1:T1}
Let $(M,\PD M, g)$ be a simple Riemannian manifold and assume that the metric $g$ is real-analytic. Let $\Ac$ be an open set of geodesics in $\wt{M}$ such that that each geodesic $\g\in \Ac$ can be deformed to a point on the boundary $\PD \wt{M}$ by geodesics in $\Ac$. Let $M_{\Ac}$ be the set of points lying on the intersection of these geodesics with $M$. If $f\in L^{2}(M)$ is a function such that $If(\g)=0$ for each geodesic $\g$ in $M$ such that its geodesic extension to $\wt{M}$ belongs to $\Ac$, then $f=0$ on $M_{\Ac}$.
\end{thm}

Note that the function is a priori defined only on $M$. We will extend $f$ to $\wt{M}$ such that $f=0$ in $\wt{M}\setminus M$. Then we will still have $If(\g)=0$ for each $\g$ in $\Ac$. 

\section{Analytic regularity of $f$ along conormal directions of $\g\in\Ac$}
We fix a geodesic $\g_{0}\in \Ac$. Since $\Ac$ is open, for each $\g$ near $\g_{0}$, we have $If(\g)=0$. Under this assumption, we show that $f$ is microlocally analytic along conormal directions of $\g_{0}$. Our references for analytic microlocal analysis are \cite[Chap. 8]{H1}, \cite[Chap.5]{T} and \cite{Sj}.
\begin{prop}\label{S2:P1}
Let $\g_{0}\in \Ac$. Then $WF_{a}(f)\cap N^{*}\g_{0}=\emptyset$.
\end{prop}
\bpr 
The proof of this proposition is fairly straightforward for $C^{\infty}$ wavefront sets. For, one can consider a $C_{c}^{\infty}$ cut-off $\psi$ in a small neighborhood of $\g_{0}$ and consider the pseudodifferential operator $I^{*}\psi I$. This is elliptic in $N^{*}\g_{0}\setminus\{0\}$. If we assume that $If=0$ near $\g_{0}$, then by microlocal elliptic regularity, we get that $f$ is smooth in $N^{*}\g_{0}\setminus\{0\}$. The proof for the analytic case is much more involved because with the cut-off $\psi$, $I^{*}\psi I$ is not an analytic pseudodifferential operator anymore and there are restrictions on the kind of cut-offs we can use. 

In the analytic case, we use a recent result of Stefanov-Uhlmann given in \cite[Prop.2]{SU3} based on the characterization of analytic singularities using the generalized FBI transform and the complex stationary phase method. See the book of J. Sj\"ostrand \cite{Sj} for more details. Actually Stefanov-Uhlmann in \cite[Prop.2]{SU3} have proved the statement of the proposition with $f$ replaced by a solenoidal symmetric tensor field. Their proof carries over here with straightforward modifications and in fact the arguments for functions are simpler. We will not repeat the proof here, but just mention that the simplification in our case is due to the fact that only one phase function satisfying the estimate in \cite[Eq. 51]{SU3} is enough to prove our result. This estimate also shows that the conormals to $\g_{0}$ are not in the analytic wavefront set of $f$; see \cite[Defn 6.1]{Sj}. We also mention a recent work of Frigyik-Stefanov-Uhlmann in \cite{FSU} where similar ideas based on the complex stationary phase method are used.
\epr
It is of interest, whether one can prove the statement of Proposition \ref{S2:P1} relying on analytic pseudodifferential operator methods alone. We are able to provide a proof for 2-dimensional manifolds. We parameterize the geodesics of $M$ by $\PD_{+}SM:=\{(x,\xi): x\in \PD M, |\xi|=1, \langle \nu,\xi\R\leq 0\}$ where $\nu$ is the outer unit normal to the boundary. Let $\g_{0}$ intersect the boundary $\PD M$ at a point $x_{0}$ in the direction $\xi_{0}$ such that  $(x_{0},\xi_{0})\in \PD_{+}S M$. We now perturb the point $x_{0}$ along the boundary and the vector $\xi_{0}$ slightly such that the geodesics starting from these perturbed points and directions still lie in $\Ac$. Let $\psi$ be a $C_{c}^{\infty}$ function supported in this perturbation and identically $1$ in a smaller neighborhood of $(x_{0},\xi_{0})$. As in \cite{BQ1}, we write
\[ I^{*}If=I^{*}\psi If+I^{*}(1-\psi)If.\]
By hypothesis, the first term on the right above is $0$. For the second term above, we analyze what $I^{*}$ does to the analytic singularities of $(1-\psi)If$. For this we make use of a calculation proved in \cite[ pp. 212, equation 2.15]{GU}. Let $Z=\{(x,y,\eta)\in M\times \Ac: \exp_{y}(s\eta)=x, (y,\eta)\in \PD SM, s\in \Rb\}$. (Here $SM$ is the unit sphere bundle). We then have that $N^{*}Z$ is given by
\begin{align}\label{S4:EQ1}
N^{*}Z=&\{(\exp_{y}(s\eta),y,\eta; \notag\\
&((D_{\eta} \exp)^{*})^{-1}(\eta\lrcorner \Omega),((D_{y}\exp)^{*}(D_{\eta}\exp)^{*})^{-1}(\eta\lrcorner \Omega),s\eta\lrcorner \Omega|_{T\Ac}):\notag\\
& (y,\eta)\in \Ac, s\in \Rb, \Omega\in \Lambda^{2}T_{y}^{*}M)\}.
\end{align}
Now we show that the projection $N^{*}Z\to T^{*}\Ac$ is injective. Suppose we are given: 
\[ (y,\eta,(D_{y}\exp)^{*}((D_{\eta}\exp)^{*})^{-1}(\eta\lrcorner \Omega), s\eta\lrcorner \Omega|_{T\Ac}). \] 
Then $s((D_{y}\exp)^{*}(D_{\eta}\exp)^{*})^{-1}(\eta\lrcorner \Omega)=((D_{y}\exp)^{*}(D_{\eta}\exp)^{*})^{-1}(s\eta\lrcorner \Omega)$. The right hand side is known and the vector part of the left hand side is known. So we can recover $s$, from which we can recover $\exp_{y}(s\eta)$.

For simplicity let us denote $u=(1-\psi)If$. Based on the injectivity of the projection shown above, we have that an analytic singularity at $\g_{0}$ pulls back to an analytic singularity at a point on the geodesic $\g_{0}$ and no other geodesic. From formula \eqref{S4:EQ1}, we get that $N^{*}Z$ does not contain any elements where one of the conormal directions is $0$. Now using \cite[Theorem 8.5.5]{H1}, we get that if the projection pulls back to an analytic singularity on the geodesic $\g_{0}$, then such a singular direction has to be conormal to the geodesic $\g_{0}$. Now $u$ is analytic at $\g_{0}$ in all directions and in $2$ dimensions the singularities of the cut-off $\psi$ pull back to singularities not conormal to the geodesic $\g_{0}$. Therefore we get that $WF_{a}(I^{*}If)\cap N^{*}\g_{0}=\emptyset$. Now by the work of Stefanov-Uhlmann in \cite{SU2}, we have that $I^{*}I$ is an elliptic analytic pseudodifferential operator. Therefore using elliptic regularity, we get $WF_{a}f\cap N^{*}\g_{0}=\emptyset$. This provides another proof of the proposition for the $2$-dimensional case. 

In higher dimensions, this proof does not work because the analytic singularities of $(1-\psi)$ could appear conormal to a different geodesic $\g_{1}$ that intersects $\g_{0}$ and it is possible that the analytic singular direction is at the point of intersection and also conormal to both $\g_{0}$ and $\g_{1}$. It is an interesting question whether by a very careful choice of the cut-off $\psi$  and some modification of the arguments above, a similar proof can be given for the higher dimensional case.

\section{Proof of the main result}
In this section, we prove the main result. Recall that we have extended the function $f$ to be $0$ in $\wt{M}\setminus M$. 
\bpr [Proof of Theorem \ref{S1:T1}]
Assume first that $n\geq 3$. Fix a geodesic $\g \in \Ac$. We define a continuous map $Q:[0,1]\to \Ac$ such that $Q(1)=\sigma$, where $\sigma\in \Ac$ is a geodesic that does not intersect $\text{supp } f$ and $Q(0)=\g$. Note that by the openness of $\Ac$, there is always a geodesic in $\Ac$ that does not intersect $\mbox{supp}(f)$.
 Now consider a geodesic $Q(t)$. Assume that this geodesic is parameterized by $(x,\xi)\in \PD S_{x}\wt{M}$ where $x$ is a point on the boundary $\PD \wt{M}$. Now consider all unit vectors that are at a fixed angle to $\xi$. If we assume this angle is small enough, then the geodesics starting at $x$ in this cone of directions lie in $\Ac$ and also that there is at least one cone that does not intersect supp $(f)$. Let us call this cone of geodesics $\Cc_{Q(t)}$. We do this construction for all $0\leq t\leq 1$. We define a new map $P$ that associates for every $0\leq t\leq 1 $, the cone $\Cc_{Q(t)}$.  Following \cite{BQ1}, we let,
\[ t_{1}=\inf\{ t: P(t_{2})\cap \text{supp } f =\emptyset \text{ for all } t_{2}>t\}.\]
 Assume $t_{1}>0$. Then $P(t_{1})$ intersects support of $f$ at the boundary on a compact set. Since the cone $P(t_{1})$ is an embedded hypersurface, by \cite[Defn 8.5.7]{H1}, we have that at a point of intersection, the normal to the cone is also normal to $\text{supp } f$. Now by \cite[Thm 8.5.6]{H1}, such a normal is in $WF_{a}(f)$. But this is a contradiction to Proposition \ref{S2:P1}. Therefore $f=0$ in a neighborhood of this point of intersection. Since the points of intersection is a compact set, by considering a finite number of such points, we can show that $f=0$ at a positive distance away from this cone. But this contradicts the infimum value $t_{1}$ above. Hence $t_{1}=0$ and so we have $f=0$ on $Q([0,1])$. For the case of $n=2$, the geodesics themselves form hypersurfaces and the same argument as above holds without the need for the cone construction. This shows that $f=0$ on the set of points lying on the deformation. The proof of the theorem is complete.
\epr
{\it Acknowledgments:} I would like to thank Gunther Uhlmann and Plamen Stefanov for the guidance and encouragement. Plamen Stefanov pointed out errors in a previous version of this paper for which I am very grateful. Thanks also to Jan Boman for the useful E-mail conversations and Eric Todd Quinto for the fruitful discussions regarding this work.


\begin{thebibliography}{10}

\bibitem{B1}
Jan Boman.
\newblock Helgason's support theorem for {R}adon transforms---a new proof and a
  generalization.
\newblock In {\em Mathematical methods in tomography (Oberwolfach, 1990)},
  volume 1497 of {\em Lecture Notes in Math.}, pages 1--5. Springer, Berlin,
  1991.

\bibitem{B2}
Jan Boman.
\newblock Holmgren's uniqueness theorem and support theorems for real analytic
  {R}adon transforms.
\newblock In {\em Geometric analysis (Philadelphia, PA, 1991)}, volume 140 of
  {\em Contemp. Math.}, pages 23--30. Amer. Math. Soc., Providence, RI, 1992.

\bibitem{BQ1}
Jan Boman and Eric~Todd Quinto.
\newblock Support theorems for real-analytic {R}adon transforms.
\newblock {\em Duke Math. J.}, 55(4):943--948, 1987.

\bibitem{BQ2}
Jan Boman and Eric~Todd Quinto.
\newblock Support theorems for {R}adon transforms on real analytic line
  complexes in three-space.
\newblock {\em Trans. Amer. Math. Soc.}, 335(2):877--890, 1993.

\bibitem{FSU}
Bela Frigyik, Plamen Stefanov, and Gunther Uhlmann.
\newblock The {X}-ray transform for a generic family of curves and weights.
\newblock {\em J. Geom. Anal.}, 18(1):89--108, 2008.

\bibitem{GQ1}
Fulton Gonzalez and Eric~Todd Quinto.
\newblock Support theorems for {R}adon transforms on higher rank symmetric
  spaces.
\newblock {\em Proc. Amer. Math. Soc.}, 122(4):1045--1052, 1994.

\bibitem{GU}
Allan Greenleaf and Gunther Uhlmann.
\newblock Nonlocal inversion formulas for the {X}-ray transform.
\newblock {\em Duke Math. J.}, 58(1):205--240, 1989.

\bibitem{GQ2}
Eric~L. Grinberg and Eric~Todd Quinto.
\newblock Morera theorems for complex manifolds.
\newblock {\em J. Funct. Anal.}, 178(1):1--22, 2000.

\bibitem{GS}
Victor Guillemin and Shlomo Sternberg.
\newblock {\em Geometric asymptotics}.
\newblock American Mathematical Society, Providence, R.I., 1977.
\newblock Mathematical Surveys, No. 14.

\bibitem{H1}
Lars H{\"o}rmander.
\newblock {\em The analysis of linear partial differential operators. {I}}.
\newblock Classics in Mathematics. Springer-Verlag, Berlin, 2003.
\newblock Distribution theory and Fourier analysis, Reprint of the second
  (1990) edition [Springer, Berlin; MR1065993 (91m:35001a)].

\bibitem{Q1}
Eric~Todd Quinto.
\newblock Real analytic {R}adon transforms on rank one symmetric spaces.
\newblock {\em Proc. Amer. Math. Soc.}, 117(1):179--186, 1993.

\bibitem{Q2}
Eric~Todd Quinto.
\newblock Support theorems for the spherical {R}adon transform on manifolds.
\newblock {\em Int. Math. Res. Not.}, pages Art. ID 67205, 17, 2006.

\bibitem{Sh1}
V.~A. Sharafutdinov.
\newblock {\em Integral geometry of tensor fields}.
\newblock Inverse and Ill-posed Problems Series. VSP, Utrecht, 1994.

\bibitem{Sj}
Johannes Sj{\"o}strand.
\newblock Singularit\'es analytiques microlocales.
\newblock In {\em Ast\'erisque, 95}, volume~95 of {\em Ast\'erisque}, pages
  1--166. Soc. Math. France, Paris, 1982.

\bibitem{SU2}
Plamen Stefanov and Gunther Uhlmann.
\newblock Boundary rigidity and stability for generic simple metrics.
\newblock {\em J. Amer. Math. Soc.}, 18(4):975--1003 (electronic), 2005.

\bibitem{SU3}
Plamen Stefanov and Gunther Uhlmann.
\newblock Integral geometry of tensor fields on a class of non-simple
  riemannian manifolds.
\newblock {\em Amer. J. Math.}, 130(1):239--268, 2008.

\bibitem{T}
Fran{\c{c}}ois Tr{\`e}ves.
\newblock {\em Introduction to pseudodifferential and {F}ourier integral
  operators. {V}ol. 1}.
\newblock Plenum Press, New York, 1980.
\newblock Pseudodifferential operators, The University Series in Mathematics.

\bibitem{QZ}
Yiying Zhou and Eric~Todd Quinto.
\newblock Two-radius support theorems for spherical {R}adon transforms on
  manifolds.
\newblock In {\em Analysis, geometry, number theory: the mathematics of Leon
  Ehrenpreis (Philadelphia, PA, 1998)}, volume 251 of {\em Contemp. Math.},
  pages 501--508. Amer. Math. Soc., Providence, RI, 2000.

\end{thebibliography}
\end{document}